\documentclass[11pt, a4paper]{amsart}
\usepackage{fullpage}
\usepackage{amsmath}
\usepackage{amssymb}
\usepackage[colorlinks,linkcolor=blue,filecolor=blue,citecolor=blue,urlcolor=blue,bookmarks=false,hypertexnames]{hyperref}
\usepackage{tikz}
\usepackage{subcaption}
\usepackage[T1]{fontenc}
\usepackage{paralist}

\usepackage[normalem]{ulem}

\newcommand\naturals{\mathbb{N}}

\def\opn#1#2{\def#1{\operatorname{#2}}} 

\let\union=\cup
\let\dirsum=\oplus
\opn\cocomp{co-comp}
\opn\comp{comp}
\opn\rank{rank}
\opn\maxind{max-ind}

\makeatletter
\@addtoreset{equation}{section}
\def\theequation{\thesection.\@arabic \c@equation}
\makeatother

\theoremstyle{plain}

\newtheorem{theorem}[equation]{Theorem}
\newtheorem{corollary}[equation]{Corollary}

\newtheorem{Lemma}[equation]{Lemma}
\newtheorem{proposition}[equation]{Proposition}

\theoremstyle{definition}

\newtheorem{definition}[equation]{Definition}

\newtheorem{discussion}[equation]{Discussion}

\newtheorem{Remark}[equation]{Remark}
\newenvironment{remarkbox}[1][]{%
\begin{Remark}[#1]\pushQED{\qed}}{\popQED \end{Remark}}

\newtheorem{example}[equation]{Example}

\title {Cohen-Macaulay permutation graphs}

\author{P.~V.~Cheri}
\address{Pondicherry University, Chinna Kalapet, Kalapet, Puducherry 605014. India}
 \email{cherivachaz@gmail.com}

\author{Deblina Dey}
\address{Indian Institute of Technology Madras, Chennai, 600036. India}
\email{deblina.math@gmail.com}

\author{Akhil K}
\address{Pondicherry University, Chinna Kalapet, Kalapet, Puducherry 605014. India}
 \email{akhilsekger@gmail.com}

\author{Nirmal Kotal}
\address{Chennai Mathematical Institute, Siruseri, Tamilnadu 603103. India}
\email{nirmal@cmi.ac.in}

\author{Dharm Veer}
\address{Indian Institute of Technology Gandhinagar, Gandhinagar, Gujarat 382355. India}
\email{dharm.v@iitgn.ac.in}

\thanks{DD was supported by the Prime Minister Research Fellowship, India. NK was partly supported by an Infosys Foundation fellowship.}



\begin{document}

\begin{abstract}
    In this article, we characterize Cohen-Macaulay permutation graphs.
    In particular, we show that a permutation graph is Cohen-Macaulay if and only if 
 it is well-covered and there exists a unique way of partitioning its vertex set into $r$ disjoint maximal cliques, where $r$ is the cardinality of a maximal independent set of the graph. 
We also provide some sufficient conditions for a comparability graph to be a uniquely partially orderable (UPO) graph.
\end{abstract}

\maketitle
\section{Introduction}

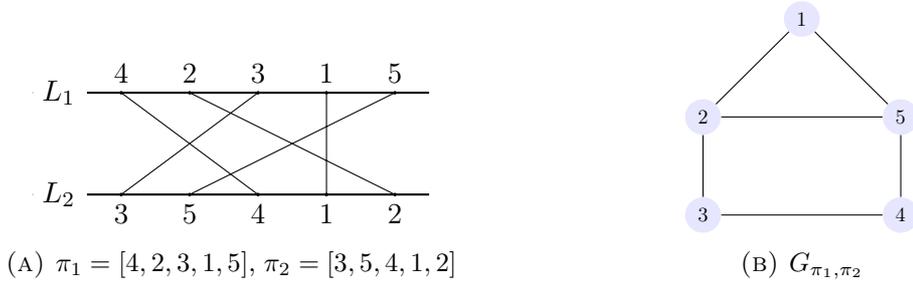
\begin{figure}[h]%
	\begin{subfigure}[t]{7cm}
		\begin{center}
		\begin{tikzpicture}[scale=.9]
			\draw[thick] (.5,.5)--(5.5,0.5)
			(.5,2)--(5.5,2);

			\draw[]      (5,2)--(2,.5)
			(2,2)--(5,0.5) 
			(1,2)--(3,0.5)
			(3,2) -- (1,0.5)
			(4,2) -- (4,0.5)
			;
			
			\filldraw[black] (-0.3,2) circle (0pt) node[anchor=west]  {$L_1$};
			\filldraw[black] (-.3,0.5) circle (0pt) node[anchor=west]  {$L_2$};
			
			\filldraw[black] (1,0.5) circle (.5pt) node[anchor=north] {3};
			\filldraw[black] (2,0.5) circle (.5pt) node[anchor=north] {5};
			\filldraw[black] (3,0.5) circle (.5pt) node[anchor=north] {4};
			\filldraw[black] (4,0.5) circle (.5pt) node[anchor=north]  {1};
			\filldraw[black] (5,0.5) circle (.5pt) node[anchor=north]  {2};
			
			\filldraw[black] (1,2) circle (.5pt) node[anchor=south] {4};
			\filldraw[black] (2,2) circle (.5pt) node[anchor=south] {2};
			\filldraw[black] (3,2) circle (.5pt) node[anchor=south] {3};
			\filldraw[black] (4,2) circle (.5pt) node[anchor=south]  {1};
			\filldraw[black] (5,2) circle (.5pt) node[anchor=south]  {5};
		\end{tikzpicture}
		\caption{\centering
			$\pi_1=[4,2,3,1,5]$,
			$\pi_2=[3,5,4,1,2]$}
		\end{center}
	\end{subfigure}
	\quad
	\begin{subfigure}[t]{7cm}
		\begin{center}
		\begin{tikzpicture}  
			[scale=1.3,auto=center,every node/.style={circle,fill=blue!10,scale=.7}]

			\node (a1) at (2,3)  {1};  
			\node (a2) at (1,2)  {2};  
			\node (a3) at (1,1)  {3}; 
			\node (a4) at (3,1)  {4};  
			\node (a5) at (3,2) {5}; 
			\draw (a1) -- (a2) --(a3) -- (a4) -- (a5)--(a1);
			\draw (a2)--(a5);  
		\end{tikzpicture}  
		\caption{$G_{\pi_1,\pi_2}$}%
			\end{center}
	\end{subfigure}

	\caption{Example of a permutation graph}
	\label{fig:PermuGraphButNotChordalBipartite}
\end{figure}

Permutation graphs were introduced in \cite{PLE71permutationgraphs,Even1972PermutationGA}. 
Let $L_1$ and $L_2$ be two horizontal lines labeled from left to right by some permutations $\pi_1$ and $\pi_2$ on $[n]=\{1,\ldots,n\}$.
For each $1\leq i \leq n,$ let $f_i$ denotes the line segment joining $i$ on $L_1$ to $i$ on $L_{2}$. A permutation graph $G_{\pi_1,\pi_2}$ is a graph on $[n]$, such that $\{i,j\}$ is an edge of $G_{\pi_1,\pi_2}$ if and only if $f_i$ intersects with $f_j$. A simple graph $G$ on $n$ vertices is a \emph{permutation graph} if there exist two permutations $\pi_1,\pi_2$ on $[n]$, such that $G=G_{\pi_1,\pi_2}$.
Figure~\ref{fig:PermuGraphButNotChordalBipartite} gives an example of a permutation graph. 

It is known that a graph is a permutation graph if and only if it is a comparability graph and a co-comparability graph~\cite[Theorem\ 3]{PLE71permutationgraphs}.
    Thus, the complement of a permutation graph is again a permutation graph. 
    It follows from~\cite{Gallai67} that a co-comparability graph does not contain induced cycles of length $n\geq 5$ (see~\cite{Gallai01translation} for an English translation of~\cite{Gallai67}).
    Therefore, a permutation graph $G$ is a weakly chordal graph, i.e., both $G$ and its complement do not have an induced cycle of length $\geq 5$.

Let $G$ be a simple graph on $[n]$.
Let $S = K[x_1,\ldots,x_n]$ be a polynomial ring over a field $K$.
Let $I_G = (x_ix_j : \{i,j\} \ \text{is an edge of} \ G)$ be the {\em edge ideal} of $G$.
We say that $G$ is {\em Cohen-Macaulay} over $K$ if $S/I_G$ is Cohen-Macaulay over $K$,
and we say $G$ is {\em Cohen-Macaulay} if $S/I_G$ is Cohen-Macaulay over any field.
Edge ideals were introduced by Villarreal in~\cite{Vil90cmgraphs}, and in the same paper, he classified all Cohen-Macaulay trees~\cite[Corollary\ 2.5]{Vil90cmgraphs}. 
Faridi~\cite{Faridi2003facetideals} extended the study of tree graphs to simplicial complexes and showed that the facet ideals of simplicial trees are unmixed if and only if Cohen-Macaulay.
Herzog and Hibi~\cite{HH05} generalized this result by characterizing Cohen-Macaulay bipartite graphs.
Crupi et al.~\cite{CRT11cmedgeideals} studied Cohen-Macaulay graphs whose height is half of the number of vertices.
 Herzog et al.~\cite{HHZ06cmchordalgraphs} proved that a chordal graph is Cohen-Macaulay if and only if it is unmixed (equivalently well-covered).  
    They also noted that it is just as challenging to describe all Cohen-Macaulay graphs as it is to describe all Cohen-Macaulay simplicial complexes. 
    Therefore, one can not anticipate a general classification theorem.
    Unmixed property of edge ideal of clutters have been studied in~\cite{MRV08clutter}; clutters are generalization of graphs.

    In this paper, we characterize Cohen-Macaulay permutation graphs and show that well-covered is not a sufficient condition for a permutation graph to be Cohen-Macaulay. 
    Thus, for a weakly chordal graph to be Cohen-Macaulay, well-covered is not a sufficient condition (see Example \ref{eg:noncm_perm}). 
    Our main theorem for Cohen-Macaulay permutation graphs is as follows.


\begin{theorem}
	\label{theorem:cm_iff_unique_partition}
    Let $G$ be a permutation graph. Then the following are equivalent:
    \begin{enumerate}
        \item
            $G$ is Cohen-Macaulay.

        \item\label{theorem:cm:unique_partition}
            $G$ is well-covered and there exists a unique way of partitioning $V(G)$ into $r$ disjoint maximal cliques, 
            where $r$ is the cardinality of a maximal independent set of $G$.
   
    \end{enumerate}
\end{theorem}

Observe that in~\eqref{theorem:cm:unique_partition} of the Theorem~\ref{theorem:cm_iff_unique_partition}, well-covered implies that the cardinality of all maximal independent sets of $G$ is equal, and it is equal to the Krull-dimension of the ring $S/I_G$. 
Proof of Theorem~\ref{theorem:cm_iff_unique_partition} is based on a result of Jahangir and the last author \cite{JV23cmposets}, that characterizes Cohen-Macaulay poset of dimension two and the fact that a permutation graph is a co-comparability graph of a poset of dimension at most two \cite{millerposetdimension}. We refer the reader to the Preliminaries section for poset related terms.

We now move towards our main result about UPO graphs.
We denote the comparability graph of a poset $P$ by $\comp(P)$.
A comparability graph $G$ is a {\em uniquely partially orderable (UPO)} graph if
$G = \comp(P) = \comp(Q)$ for some posets $P$ and $Q$ implies $P = Q$ or $P = Q^{\partial}$ where $Q^\partial$ denotes the dual poset of $Q$.
M\"{o}hring~\cite{moh84upo} proved that almost all comparability graphs are UPO, i.e., $lim_{n\to \infty} G(n, UPO)/G(n)=1$, where $G(n)$ and $G(n, UPO)$ denote the numbers of comparability
graphs and UPO on the vertex set $[n]$.
In this direction, we prove the following result.

 \begin{theorem}
     \label{thm:upo}
     Let $P$ be a pure poset such that the induced subposet of $P$ consisting of height $i$ and height $j$ elements is connected for all $0 \leq i<j \leq \rank(P)$.  
     If $P$ is not an ordinal sum of two induced subposets, then $\comp(P)$ is a UPO.
 \end{theorem}

 To prove Theorem~\ref{thm:upo}, we use a result of Trotter et~al.~\cite[Theorem\ 1]{TMS76dimensioncomparability} on UPO which shows that 
 a connected comparability graph is a UPO if and only if every nontrivial partitive subset is an independent set of vertices. 
 We refer the reader to Section~\ref{sec:thm:upo} for the definition of a partitive subset. 
 As a corollary of Theorem~\ref{thm:upo}, we show that the complement of a Cohen-Macaulay connected permutation graph is a UPO.

 Section~\ref{sec:pre} contains the definitions and preliminaries. The proof of the Theorem~\ref{theorem:cm_iff_unique_partition} is given in Section~\ref{sec:thm:cm}. 
Section~\ref{sec:thm:upo} is devoted to the proof of Theorem~\ref{thm:upo} 

\section{Preliminaries}\label{sec:pre}

Let $(P,\leq)$ be a poset. We say that two elements $x$ and $y$ of $P$ are {\em comparable} if $x\leq y$ or $y \leq x$; otherwise, $x$ and $y$ are {\em incomparable}. For $x,y\in P$, we say that $y$ {\em covers} $x$, denoted by $x \lessdot y$, if $x<y$ and there is no $z\in P$ with $x<z<y$. 
We say that $P$ is an {\em antichain} if any two distinct elements of $P$ are incomparable. 
The {\em dual poset} of $(P,\leq)$, denoted by $(P^\partial,\preceq)$, is the poset with the same underlying set, but its order relations are opposite of $P$, i.e. $p\leq q$ if and only if $q\preceq p.$ 
The {\em Hasse diagram} of a poset $P$ is the graph whose vertices are elements of $P$, whose edges are cover relations and such that if $x<y$ the $y$ is "above" $x$ (i.e., with a higher vertical coordinate). 
 
A {\em chain} $C$ of $P$ is a totally ordered subset of $P$. 
 The length of a  chain $C$ of $P$ is $\#C-1$, where $\#C$ is the cardinality of $C$. The {\em rank} of $P$, denoted by rank$(P)$, is the maximum of the lengths of chains in $P$. 
 A poset is called {\em pure} if all maximal chains of $P$  have the same length. 
 An induced subposet $Q$ of $P$ is a poset on a subset of the underlying set $P$ such that for every $x,y\in Q,$ $x\leq y$ in $Q$ if and only if $x\leq y$ in $P$. 
 For $p\in P$, the {\em height} of $p$ is the rank of the induced subposet of $P$ which consists of  all $q \in P$ with $q \leq p$.

Let $G$ be a simple graph with vertex set $V(G)$ and edge set $E(G)$. 
For any $S\subseteq V(G)$, $G[S]$ denotes the induced subgraph of $G$ on the set $S$. 
A {\em clique} of $G$ is a subset $S$ of the vertex set of $G$ such that for any $i,j \in S$, $\{i,j\} \in E(G)$. 
A {\em maximal clique} is a clique of $G$ of maximum size, i.e., it is not contained in any other clique with respect to inclusion. 
A vertex \( v \in V(G) \) is called a \emph{simplicial vertex} if its neighborhood \( G\left[\mathcal{N}(v)\right] \), where
	\(
	\mathcal{N}(v) := \{u \in V(G) \mid \{u,v\} \in E(G)\},
	\)
	forms a clique.
A set $S \subseteq V(G)$ is an {\em independent set} of $G$ if for any $i,j\in S$, $\{i,j\} \notin E(G)$. 
An independent set of maximal size is called a maximal independent set. 
The graph $G$ is {\em well-covered} if all maximal independent sets of $G$ are of the same size. 
The cardinality of a maximal independent set of a graph $G$ is known to be the Krull-dimension of $S/I_G$.

A {\em co-comparability graph} (resp.~{\em comparability graph}) $G$ corresponding to a poset $P$ is the graph with the vertex set same as the underlying set of $P$, and $\{x,y\} \in E(G)$ if and only if $x$ and $y$ are incomparable (resp.~comparable and $x\neq y$) in $P$. 
A graph is said to be co-comparability (resp.~comparability graph) if it is a co-comparability graph (resp.~comparability graph) of some poset.
We denote the co-comparability graph (resp.~comparability graph) of a poset $P$ by $\cocomp(P)$ (resp.~$\comp(P)$).
Note that the comparability graph of a poset $P$ is the complement of the co-comparability graph of $P$.

The order complex of a poset $P$ is a simplicial complex on the underlying set of $P$ whose faces are chains of $P$.   
The poset $P$ is said to be {\em Cohen-Macaulay} over a field $K$ if the Stanley-Riesner ring associated with its order complex is Cohen-Macaulay over $K$.
It is known that the Stanley-Reisner ideal of the order complex of $P$ coincides with the edge ideal of the $\cocomp(P)$.
Thus, $P$ is Cohen-Macaulay over a field $K$ if and only if $\cocomp(P)$ is Cohen-Macaulay over $K$.

A linear extension $\pi$ of a poset  $P$ is a linear order on the underlying set of $P$ such that $x\leq y$ in $\pi$ whenever $x\leq y$ in $P$. A  poset $P$ is an intersection of a family of linear extensions $\pi_{1},\ldots ,\pi_{d}$ if $x\leq y$ in $P$ if and only  if $x\leq y$ in $\pi_{i}$, for all $1\leq i \leq d$. The {\em dimension of a poset} $P$ is the least integer $d$ such that $P$  can be expressed as the intersection of $d$ linear extensions of $P$. The dimension of a poset was defined by Dushnik and Miller~\cite{millerposetdimension}.   
The Cohen-Macaulay posets of dimension two were characterized in~\cite[Theorem\ 1]{JV23cmposets}. 

Let $P$ and $Q$ be two posets on disjoint sets. The disjoint union of $P$ and $Q$ is the  poset $P+Q$ on the set $P\cup Q$ with the following order: if $x,y\in P+Q,$ then $x\leq y$ if either  $x,y\in P$ and $x \leq y$ in $P$ or $x,y \in Q$ and $x\leq y $ in $Q$. A poset which can be written as a disjoint union of two posets is called {\em disconnected}; otherwise, the poset is called {\em connected}.
 The {\em ordinal sum} $P\dirsum Q$ of $P$ and $Q$ is the poset on the set $P\cup Q$ with the following order: if $x, y \in P\dirsum Q $, then $x \leq y$ if either $x, y \in P$ and $x \leq y$ in $P$ or $x, y \in Q$ and $x \leq y$ in $Q$ or $x \in P$ and $y \in Q$. 

\section{Proof of Theorem~\ref{theorem:cm_iff_unique_partition}}\label{sec:thm:cm}

To prove Theorem~\ref{theorem:cm_iff_unique_partition}, we use a result of~\cite{JV23cmposets} that characterizes Cohen-Macaulay posets of dimension two (see Lemma~\ref{thm:cmpermutation}), and the fact that a permutation graph is a co-comparability graph of a poset of dimension at most two \cite{millerposetdimension}. 
 We start by stating a few lemmas concerning dimension two pure posets.

     We first recall the linear order defined on the height $i$ elements in~\cite{JV23cmposets} for a pure poset of dimension two.
    Let $P$ be a pure poset on $[n]$ that is the intersection of the identity permutation and another permutation.
    For $0\leq i\leq \rank(P)$, we denote the set of height $i$ elements of $P$ by $P_i$, and
    define a linear order $<_i$ on $P_i$ as follows: for any $x,y\in P_i$, 
	\[ x<_i y ~\text{if and only if}~ x>y ~\text{in}~ \naturals. \]
For any $x \in P_i$, define $U(x)=\{y \in P: x \lessdot y\} $. 
Since $P$ is pure, $U(x)\subseteq P_{i+1}$. 
We denote the minimal (resp. maximal) element of on $U(x)$ under the linear order $P_{i+1}$ by $x_{min}$ (resp. $x_{max}$).

We need the following observation in the proof of our main result of this section. 
	
\begin{Lemma}\label{lem:dim2coverrelations}
		Let $P$ be a pure poset of dimension two on the underlying set $[n]$. Assume that $P$ is the intersection of the identity permutation and another permutation, and  
		the induced subposet of $P$ consisting of height $i$ and height $i+1$ elements is connected for all $0 \leq i \leq \rank(P) -1$.
		If $x \lessdot_i y$ in $P_i$ for some $x,y\in P_i$ and $0\leq i\leq \rank(P)-1$, then $y_{\min} \leq_{i+1} x_{\max}$ in $P_{i+1}$.
\end{Lemma}
\begin{proof}
	Proved as a claim in the proof of~\cite[Lemma\ 3.6]{JV23cmposets}.
\end{proof}

The following lemma states the Cohen-Macaulay characterization of posets of dimension two~\cite[Theorem\ 1]{JV23cmposets}. 
We also add an equivalent condition, which will be needed in the last section.

\begin{Lemma}\label{thm:cmpermutation}
	Let $P$ be a finite poset of dimension two.
	Then the following are equivalent:
	\begin{enumerate}
		
		\item \label{thm:cmpermutation:cm}
		$P$ is Cohen-Macaulay over any field $K$.
		
		\item \label{thm:cmpermutation:poset}
		$P$ is an antichain; or $P$ is pure and the induced subposet of $P$ consisting of height $i$ and height $i+1$ elements is connected for all $0 \leq i \leq rank(P) -1$.  
		
		\item \label{thm:cmpermutation:poset_all}
		$P$ is an antichain; or $P$ is pure and the induced subposet of $P$ consisting of height $i$ and height $j$ elements is connected for all $0 \leq i<j \leq rank(P)$.  
		
	\end{enumerate}
\end{Lemma}
\begin{proof}
	\eqref{thm:cmpermutation:cm} if and only if \eqref{thm:cmpermutation:poset} was proved in \cite[Theorem\ 1]{JV23cmposets}. 
	\eqref{thm:cmpermutation:poset_all} $\implies$ \eqref{thm:cmpermutation:poset} is obvious. 
	
	For \eqref{thm:cmpermutation:poset} $\implies$ \eqref{thm:cmpermutation:poset_all}, assume that \eqref{thm:cmpermutation:poset} holds.
	For some $0 \leq i<j \leq \rank(P)$, let $P_{i,j}$ be the induced subposet of $P$ consisting of height $i$ and height $j$ elements.
	Since $P$ is pure, for every $x\in P_i$, there exists a $y\in P_j$ such that $x< y$ in $P$ and for every $y\in P_j$ there exists an $x\in P_i$ such that $x< y$ in $P$; thus $P_{i,j}$ is pure.   
    By~\cite[Proposition\ 3.2]{JV23cmposets}, we may assume that $P$ is the intersection of the identity permutation and another permutation. 

    Since $P_{i,j}$ is pure, it suffices to show that any two elements of $P_i$ are connected in the Hasse diagram. 
    Thus, it is enough to show that if $x \lessdot_i y$ in $P_i$, then there exists a $z\in P_j$ such that $x<z$ and $y<z$ in $P_{i,j}$.
	Let $x,y\in P_i$ be with $x \lessdot_i y$ in $P_i$. 
	Then, by Lemma~\ref{lem:dim2coverrelations}, there exists a $z'\in P_{i+1}$ such that $x<z'$ and $y<z'$ in $P$. 
	Since $P$ is pure, there exists a $z\in P_j$ such that $z'\leq z$ in $P$. Thus, $x<z$ and $y<z$ in $P$; hence in $P_{i,j}$.
	This completes the proof. 
\end{proof}

The following proposition, which describes the co-comparability graph of a pure poset, might be known, but we could not find a reference. We provide a proof here as we will need it later.

\begin{proposition}
\label{thm:wellcovered}
Let $G$ be a well-covered graph.
Then the following are equivalent:
\begin{enumerate}
    \item 
        $G$ is a co-comparability graph of a pure poset of rank $r$. 

    \item
        $(i)$ There exists a partition $Y_0,Y_1,\ldots,Y_r$ of $V(G)$ such that $Y_i$ is a maximal clique of $G$ for all $0\leq i \leq r$. \\
        $(ii)$ Upto a relabeling of $Y_0,\ldots,Y_r$, $\{x,y\}\notin E(G)$ for some $x\in Y_i$ and $y\in Y_j$ with $i<j$ if and only if 
        there exists a sequence $x=x_i,x_{i+1}, \ldots,x_j =y$ such that $x_k\in Y_k$ and $\{x_k,x_{k+1}\}\notin E(G)$ for all $i\leq k\leq j-1$.
\end{enumerate}
\end{proposition}

\begin{proof}
   First, assume that $G$ is a co-comparability graph of a pure poset $P$ of rank $r$. 
   For $0\leq i\leq r$, define $P_i$ to be the set of all height $i$ elements of $P$ and $Y_i$ to be the  set of vertices in $G$ corresponding to $P_i$. 
   Since $P_i$ is an antichain, $Y_i$ is a clique of $G$ for all $i$.
   Fix an $i$ with $0\leq i\leq r$. For any $x\notin Y_i$, there exists an element $y\in Y_i$ such that $x$ and $y$ are comparable in $P$ 
   because $P$ is a pure poset.
   Thus, $Y_i\cup \{x\}$ is not a clique of $G$ for any $x\in V(G)\setminus Y_i$; thus, $Y_i$ is a maximal clique of $G$.
   Therefore, $V(G)$ is a disjoint union of maximal cliques $Y_0,\ldots,Y_r$.

   Now, assume that $\{x,y\}\notin E(G)$ for some $x\in Y_i$ and $y\in Y_j$ with $i<j$.
   Then, $x$ and $y$ are comparable in $P$. 
   Since $P$ is pure, there exists a chain $x=x_i \lessdot x_{i+1} \lessdot \cdots \lessdot x_j=y$ in $P$ with $x_k\in P_k$ for all $i\leq k\leq j$.
   Thus, $\{x_k,x_{k+1}\}\notin E(G)$ for all $i\leq k\leq j-1$.
   On the other hand, suppose that there exists a sequence $x=x_i, x_{i+1}, \ldots, x_j=y$ such that $x_k\in Y_k$ and $\{x_k,x_{k+1}\}\notin E(G)$, $i \leq k \leq j-1$. 
   Then, $x_k$ and $x_{k+1}$ are comparable in $P$ for each $k$. Thus, $x$ and $y$ are comparable in $P$. So, $\{x,y\} \notin E(G)$.

   For the converse, we show that there exists a pure poset $P$ on the underlying set $V(G)$ such that $\cocomp(P) = G$.
   Let $Y_{0},Y_{1},\ldots,Y_{r}$ be the labeling under which $(ii)$ holds.
   We define the poset $P$ on $V(G)$ with the following cover relations: 
   for $x,y\in V(G)$, $x\lessdot y$ in $P$ if $x\in Y_{i}, y \in Y_{i+1}$ and $\{x,y\} \notin E(G)$ for some $0\leq i \leq r-1$.
   Observe that the transitivity of $P$ holds from $(ii)$.
   Clearly, $Y_i$ is an antichain in $P$ for all $0\leq i \leq r$.

   We show that $G^c = \comp(P)$ and $P$ is pure to complete the proof. 
   Let $\{x,y\} \in E(G^c)$ with $x\in Y_i$ and $y\in Y_j$ with $i<j$. 
   Then by $(ii)$, there exists a sequence $x=x_i,x_{i+1}, \ldots,x_j =y$ such that $x_k\in Y_k$ and $\{x_k,x_{k+1}\}\in E(G^c)$ for all $i\leq k\leq j-1$.
   Note that by definition of $P$, $x_k\lessdot x_{k+1}$ in $P$ for all $i\leq k\leq j-1$. Hence, $x\leq y$ in $P$.

   On the other hand, let $x\leq y$ in $P$ where $x\in Y_i$ and $y \in Y_j$ with $i <j$. 
   Then, by definition of $P$, there exists a chain $x=x_i\lessdot x_{i+1}\lessdot \cdots \lessdot x_j =y$ such that $x_k\in Y_k$ and
   $\{x_k,x_{k+1}\}\in E(G^c)$ for all $i\leq k\leq j-1$. 
   Therefore, $\{x,y\} \in E(G^c)$ by $(ii)$.
   Hence, $G^c = \comp(P)$.

   The maximal chains of the poset $P$ correspond to the maximal independent sets of $G$.
   Since $G$ is well-covered, the cardinality of maximal independent sets of $G$ are of the same cardinality. 
   Thus, the maximal chains of $P$ have the same length.
   Hence, $P$ is pure.
\end{proof}

\begin{discussion}
	\label{discussion:notation}
    Let $G$ be a well-covered co-comparability graph with cardinality of a maximal independent set is $r$. 
    By Proposition~\ref{thm:wellcovered}, $G=\cocomp(P)$ for some pure poset $P$ of rank $r-1$.
    For each $0\leq i \leq r-1$, $P_i$ denotes the set of all height $i$ elements in $P$ and $Y_i$ is the set of vertices in $G$ corresponding to $P_i$. 
    Again, by Proposition~\ref{thm:wellcovered}, each $Y_i$ is a maximal clique of $G$, and they form a partition of $V(G)$.
\end{discussion}

We are now ready to prove Theorem~\ref{theorem:cm_iff_unique_partition}.

\begin{proof}[Proof of Theorem~\ref{theorem:cm_iff_unique_partition}]
    Since $G$ is a permutation graph, $G=\cocomp(P)$ for some poset $P$ of dimension at most $2$~\cite{GRU83comparabilitygraphs}.
    Since a Cohen-Macaulay graph is well-covered, $P$ is pure under both statements. 
    It follows from the definition of the dimension of a poset and of co-comparability graph that $\dim P=1$ if and only if $P$ is a chain if and only if $G$ is a disjoint union of vertices. 
    Also, note that $\rank(P) = 0$ if and only if $P$ is an antichain if and only if $G$ is a complete graph.  
    In both cases, the proof follows immediately. So, we may assume that $\dim P=2$ and $\rank(P)\geq 1$. 
    Upto a relabeling, we may also assume that $P$ is the intersection of the identity permutation and another permutation (see~\cite[Proposition\ 3.2]{JV23cmposets}).
    Assume the notations as in Discussion~\ref{discussion:notation}. 
Thus, $V(G)$ can always be partitioned into maximal cliques $Y_0,\ldots,Y_{r-1}$.
So it is enough to show that $G$ is Cohen-Macaulay if and only if $Y_0,\ldots,Y_{r-1}$ is the only way of partitioning $V(G)$ into maximal cliques.
	
First, assume that $G$ is Cohen-Macaulay. Suppose that there exists another partition $X_0,X_1,\newline\ldots,X_{r-1}$ of $V(G)$ such that $X_i$ is a maximal clique in $G$ for all $0\leq i\leq r-1$. 
Notice that if $Y_j\subseteq X_i$ for some $i,j$, then $Y_j=X_i$ because $Y_j$ is a maximal clique. 
    Let $k$ be the least integer such that, up to a relabelling of $X_i$'s, $Y_i=X_i$ for all $i<k$ and $Y_k\neq X_j$ for any $j\geq k$. 

    In the poset $P$, if $P_k$ is a singleton, then $Y_k\subseteq X_j$ for some $j\geq k$; thus $Y_k = X_j$. Therefore $P_k$ can not be a singleton. 
    Assume that $P_k=\{a_1, a_2, \ldots, a_{s}\}$ with $a_1\lessdot_k a_2\lessdot_k \cdots \lessdot_k a_{s}$. 
    Then, there exists an integer $l_k$, such that the vertices $ a_1, a_2,\ldots ,  a_{l_k}$ belong to $X_j$ and $a_{l_k+1}$ does not belong to $X_{j}$ for some $j\geq k$. 
    Without loss of generality, we may assume that the vertices $a_1,\ldots, a_{l_k}\in X_k$ and $a_{l_k+1}\in X_{k+1}$.

    As $G$ is Cohen-Macaulay, we get that $P$ is Cohen-Macaulay, and consequently, by Lemma~\ref{thm:cmpermutation}, the induced subposet of $P$ on the underlying set $P_k\union P_{k+1}$ is connected. 
    Thus, by  Lemma~\ref{lem:dim2coverrelations}, there exists a $c_{k+1}\in P_{k+1}$ such that $c_{k+1}$ is comparable to $a_{l_k}$ and $a_{l_k+1}$. 
    Since $P$ is pure, there exists a chain $\{c_{k+1},\ldots,c_{r-1} \}$ such that $c_j\in P_j$ for all $ k+1 \leq j\leq r-1$. 
    So, $c_i$ and $c_j$ are not adjacent in $G$ for any $k+1\leq i<j \leq r-1$. 
     Also, $c_j$ is not adjacent with $a_{l_k}$ and $a_{l_{k+1}}$ in $G$ for any $j\geq k+1$. 
     Therefore, $c_{k+1},\ldots,c_{r-1}$ must belong to $X_{k+2}\sqcup X_{k+3}\sqcup \ldots \sqcup X_{r-1}$. 
     By pigeonhole principle, there exist some $i$ and $j$ with $k+1\leq i<j\leq r-1$, such that $c_i, c_j \in X_s$ for some $k+2\leq s \leq r-1$. 
     This is a contradiction because $c_i$ and $c_j$ are not adjacent in $G$. 
     Hence, the forward direction follows.

     For the converse part, suppose that $G$ is not Cohen-Macaulay. Therefore, $P$ is not Cohen-Macaulay; thus by Lemma~\ref{thm:cmpermutation}, there exists an $i$ with $0\leq i \leq r-2$ such that the induced subposet of $P$ on the underlying set $P_i\union P_{i+1}$ is not connected.
     Therefore, there are non-empty disjoint partitions of $P_i$ and $P_{i+1}$, namely
     $P_i=P_i^1\union P_i^2$ and $P_{i+1}=P_{i+1}^1\union P_{i+1}^2$ such that every element in $P_i^1$ (respectively $P_{i}^2$) is incomparable with every element of $P_{i+1}^2$ (respectively $P_{i+1}^1$). Let $Q_i=P_i^1\union P_{i+1}^2$ and $Q_{i+1}=P_{i}^2\union P_{i+1}^1$. Since $Q_i$ (respectively $Q_{i+1}$) is an antichain in $P$, corresponding vertices in $G$ forms a clique, say $X_i$ (respectively $X_{i+1}$). 

    We show that $X_i$ and $X_{i+1}$ are maximal cliques of $G$.
    If $X_i$ is not a maximal clique, then there exists a vertex $x\in V(G)\setminus X_i$ such that $ G[X_i\union x]$ is a clique. 
    If $x\in X_{i+1}$, then $x$ is adjacent to every vertex of $X_i\union X_{i+1}$ and hence $x$ is incomparable with every element of $P_i \union P_{i+1}$ in the poset $P$, which contradicts the purity of $P$. 
Thus, $x\not\in X_i\union X_{i+1}$, that is $x\not \in P_i\union P_{i+1}$ in $P$. 
As $P$ is pure, there exists a $y\in P_i$ and a $z\in P_{i+1}$ such that $\{x,y,z\}$ is a chain in $P$ in some order. 
Since $x$ is incomparable with every element of $Q_i$, $y$ and $z$ must belong to $Q_{i+1}$. 
Which is a contradiction as $Q_{i+1}$ is an antichain.
Hence, $X_i$ is a maximal clique in $G$. Similarly, $X_{i+1}$ is also a maximal clique in $G$. 
Therefore, $Y_0,\ldots,Y_{i-1},X_i,X_{i+1},Y_{i+1},\ldots,Y_{r-1}$ forms a partition of $V(G)$ into maximal cliques. This partition is different from $Y_0,\ldots,Y_{r-1}$. 
Hence, the converse direction follows by contradiction.
\end{proof}

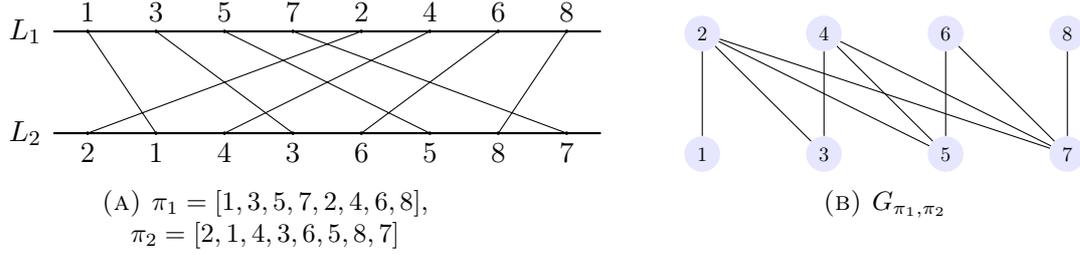
\begin{figure}[h]%
	\begin{subfigure}[t]{7cm}
		\centering
		\begin{tikzpicture}[scale=.9]
			\draw[thick] (.5,.5)--(8.5,0.5)
			(.5,2)--(8.5,2);
			
			\draw[]      (5,2)--(1,.5)
			(6,2)--(3,0.5)
			(7,2)--(5,0.5)
			(2,2)--(4,0.5) 
			(1,2)--(2,0.5)
			(3,2) -- (6,0.5)
			(8,2) -- (7,0.5)
			(4,2) -- (8,0.5)
			;
			
			\filldraw[black] (-0.3,2) circle (0pt) node[anchor=west]  {$L_1$};
			\filldraw[black] (-.3,0.5) circle (0pt) node[anchor=west]  {$L_2$};
			
			\filldraw[black] (1,0.5) circle (.5pt) node[anchor=north] {2};
			\filldraw[black] (2,0.5) circle (.5pt) node[anchor=north] {1};
			\filldraw[black] (3,0.5) circle (.5pt) node[anchor=north] {4};
			\filldraw[black] (4,0.5) circle (.5pt) node[anchor=north]  {3};
			\filldraw[black] (5,0.5) circle (.5pt) node[anchor=north]  {6};
			\filldraw[black] (6,0.5) circle (.5pt) node[anchor=north]  {5};
			\filldraw[black] (7,0.5) circle (.5pt) node[anchor=north]  {8};
			\filldraw[black] (8,0.5) circle (.5pt) node[anchor=north]  {7};
			
			\filldraw[black] (1,2) circle (.5pt) node[anchor=south] {1};
			\filldraw[black] (2,2) circle (.5pt) node[anchor=south] {3};
			\filldraw[black] (3,2) circle (.5pt) node[anchor=south] {5};
			\filldraw[black] (4,2) circle (.5pt) node[anchor=south]  {7};
			\filldraw[black] (5,2) circle (.5pt) node[anchor=south]  {2};
			\filldraw[black] (6,2) circle (.5pt) node[anchor=south]  {4};
			\filldraw[black] (7,2) circle (.5pt) node[anchor=south]  {6};
			\filldraw[black] (8,2) circle (.5pt) node[anchor=south]  {8};
		\end{tikzpicture}
		\caption{\centering
			$\pi_1=[1,3,5,7,2,4,6,8]$,
			$\pi_2=[2,1,4,3,6,5,8,7]$}
		\label{figure:permutation}
	\end{subfigure}
	\quad
    \begin{subfigure}[t]{8.3cm}
		\centering
		\begin{tikzpicture}  
            [scale=1.6,auto=center,every node/.style={circle,fill=blue!10,scale=.7}]
			\node (a8) at (4,2) {8};  
			\node (a7) at (4,1)  {7};   
			\node (a6) at (3,2)  {6};  
			\node (a5) at (3,1) {5};  
			\node (a1) at (1,1)  {1};  
			\node (a2) at (1,2)  {2};  
			\node (a4) at (2,2)  {4};  
			\node (a3) at (2,1)  {3}; 
			\draw (a1) -- (a2);  
			\draw (a2) -- (a3);  
			\draw (a2) -- (a5);  
			\draw (a4) -- (a5);  
			\draw (a3) -- (a4);  
			\draw (a6) -- (a7);  
			\draw (a5) -- (a6);  
			\draw (a7) -- (a8);
			\draw (a7) -- (a4);
			\draw (a2) -- (a7);
		\end{tikzpicture}  
		\caption{$G_{\pi_1,\pi_2}$}%
	\end{subfigure}
	\caption{A Cohen-Macaulay permutation graph}
	\label{fig:permugraph}
\end{figure}

    The graphs shown in Figures~\ref{fig:PermuGraphButNotChordalBipartite} and~\ref{fig:permugraph} are examples of Cohen-Macaulay permutation graphs.
    For better understanding of Theorem \ref{theorem:cm_iff_unique_partition}, we give an example. 

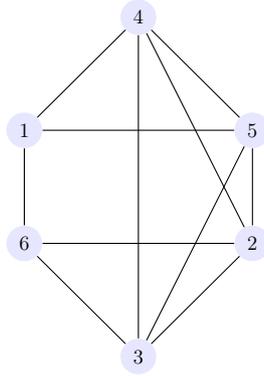
\begin{figure}[h]
	\begin{tikzpicture}  
		[scale=1.5,auto=center,every node/.style={circle,fill=blue!10,scale=.7}]
		\node (a6) at (0,1) {6};  
		\node (a2) at (2,1)  {2};   
		\node (a3) at (1,0)  {3};  
		\node (a1) at (0,2) {1};  
		\node (a5) at (2,2)  {5};  
		\node (a4) at (1,3)  {4};  
		\draw (a1) -- (a4);  
		\draw (a1) -- (a6);  
		\draw (a1) -- (a5);  
		\draw (a2) -- (a5);  
		\draw (a2) -- (a3);  
		\draw (a2) -- (a4);  
		\draw (a2) -- (a6);  
		\draw (a3) -- (a5);  
		\draw (a3) -- (a6);
		\draw (a3) -- (a4);
		\draw (a4) -- (a5);  
	\end{tikzpicture}  
\caption{A non Cohen-Macaulay permutation graph (Example~\ref{eg:noncm_perm}).}\label{figure:noncm_perm_graph}
\end{figure}

\begin{example}
	\label{eg:noncm_perm}
    Let $\pi = [5,4,6,1,3,2]$ be a permutation and $\mathrm{Id}$ be the identity permutation on six elements. 
    The permutation graph $G = G_{\mathrm{Id},\pi}$ is as shown in Figure~\ref{figure:noncm_perm_graph}. 
    Notice that $G$ is well-covered and $Y_0=\{1,4,5\},Y_1=\{2,3,6\}$ are disjoint maximal cliques of $G$. But $X_0=\{1,6\},X_1=\{2,3,4,5\}$ are also disjoint maximal cliques of $G$. 
Therefore, by Theorem~\ref{theorem:cm_iff_unique_partition}, $G$ is not Cohen-Macaulay. 
\end{example}

\begin{corollary} 
	\label{cor:zero_simplicial_imply_cm}
        Let $G$ be a well-covered permutation graph with cardinality of a maximal independent set $r$. 
        Let $Y_0,\ldots,Y_{r-1}$ be a disjoint partition of $V(G)$ into maximal cliques.
    If the number of $Y_i$'s with no simplicial vertex is at most one in each connected component of $G$, then $G$ is Cohen-Macaulay.
\end{corollary}
\begin{proof}
    It suffices to show that every connected component of $G$ is Cohen-Macaulay~\cite[Lemma\ 4.1]{Vil90cmgraphs}.
    Since an induced subgraph of a permutation graph is a permutation graph, we may assume that $G$ is connected. Suppose that $G$ is not Cohen-Macaulay. 
    From the second part of the proof of Theorem \ref{theorem:cm_iff_unique_partition}, there exist two disjoint maximal cliques $X_i$ and $X_{i+1}$ of $G$ different from $Y_i$ and $Y_{i+1}$ such that $G[X_i\union X_{i+1}]=G[Y_i\union Y_{i+1}]$ for some $0\leq i\leq r-2$.     
    As a result, each vertex of $Y_i$ belongs to two different maximal cliques, namely $Y_i$ and one of $X_i$ or $X_{i+1}$. Thus, there is no simplicial vertex in $Y_i$. 
    Similarly, $Y_{i+1}$ has no simplicial vertex. 
    In other words, there are at least two disjoint maximal cliques with no simplicial vertex, which is a contradiction.
    This completes the proof.
\end{proof}

\begin{remarkbox}
    The converse of Corollary~\ref{cor:zero_simplicial_imply_cm} may not be true. For example, take the permutation graph $G_{\pi_{1},\pi_2}$ as in Figure~\ref{fig:permugraph}. 
    Then, $G_{\pi_1,\pi_2}$ is a well-covered graph with the cardinality of a maximal independent set $4$.
    For $1\leq i \leq 4$, let $Y_i=\{2i-1,2i\}\subseteq V(G_{\pi_1,\pi_2})$. 
    Then $\{Y_1,\ldots,Y_4\}$ form a disjoint partition of $V(G_{\pi_1,\pi_2})$ into maximal cliques and such a partition is unique. 
    Hence, $G_{\pi_1,\pi_2}$ is Cohen-Macaulay by Theorem \ref{theorem:cm_iff_unique_partition}, but $Y_2$ and $Y_3$ have no simplicial vertex.
\end{remarkbox}

\begin{figure}[h]
	\begin{tikzpicture}  
		[scale=1.7,auto=center,every node/.style={circle,fill=blue!10,scale=.7}]
		\node (a1) at (0,1) {1};  
		\node (a2) at (1,1) {2};   
		\node (a3) at (2,1) {3};  
		\node (a4) at (3,1) {4};  
		\node (a5) at (4,1) {5};  
		\node (a6) at (0,0) {6};  
		\node (a7) at (1,0) {7};  
		\node (a8) at (2,0) {8};  
		\node (a9) at (3,0) {9};  
		\node (a10) at (4,0){10};  
		\draw (a1) -- (a6);  
		\draw (a1) -- (a8);  
		\draw (a1) -- (a9);  
		\draw (a1) -- (a10);  
		\draw (a2) -- (a7);  
		\draw (a2) -- (a10);  
		\draw (a3) -- (a8);  
		\draw (a3) -- (a9);  
		\draw (a3) -- (a10);
        \draw (a4) -- (a9);
		\draw (a5) -- (a10);  
	\end{tikzpicture}  
\caption{A Cohen-Macaulay bipartite non-permutation graph.}\label{figure:cm_bipartite_graph}
\end{figure}
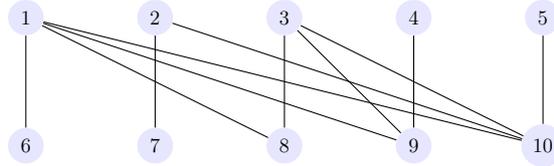

It is known that every bipartite graph is a comparability graph~\cite[Page\ 133]{BB2021linegraph}. 
By \cite[Theorem\ 3.4]{HH05}, one can deduce that a Cohen-Macaulay bipartite graph does not contain an induced cycle of length $\geq 5$.
So Cohen-Macaulay bipartite graphs are weakly chordal.
Since a permutation graph is also weakly chordal, it is natural to ask whether a Cohen-Macaulay bipartite graph a permutation graph.
The answer to this question is negative.
Consider the graph $G$ as shown in Figure~\ref{figure:cm_bipartite_graph}.
By \cite[Theorem\ 3.4]{HH05}, $G$ is Cohen-Macaulay but \texttt{SageMath}~\cite{sage} computation shows that $G$ is not a co-comparability graph.
Thus, $G$ is not a permutation graph.

\section{Proof of Theorem~\ref{thm:upo}}\label{sec:thm:upo}

Let $G$ be a graph and $K$ be a subset of $V(G)$. $K$ is said to be {\em partitive} if
for every $x\in V(G)$ with $x \notin K$, if there exists a vertex $y \in K$ such that $x$ and
 $y$ are adjacent, then $x$ is adjacent to every vertex in $K$. 
 A partitive subset $K$ is said to be nontrivial when $K$ is not the empty set, a singleton, 
 or the entire vertex set.

 Trotter et al.~\cite[Theorem\ 1]{TMS76dimensioncomparability} proved that a connected comparability graph is UPO if and only if every nontrivial
 partitive subset is an independent set of vertices.
 We use this result to prove Theorem~\ref{thm:upo}.

 \begin{proof}[Proof of Theorem~\ref{thm:upo}]
     If $\rank(P)=0$, then $P$ is an antichain; thus, $\comp(P)$ is a disjoint union of vertices which is a UPO. 
     So, we may assume that $\rank(P)=r\geq 1$.
     Observe that under the given hypothesis, $P$ is connected; thus, $\comp(P)$ is connected.
     Let $K$ be a non-empty proper subset of $V(\comp(P))$ such that it is not an independent set.
    By~\cite[Theorem\ 1]{TMS76dimensioncomparability}, it is enough to show that $K$ is not partitive.
    On the contrary, suppose that $K$ is partitive.
     Since $K$ is not independent, there exist $a, b\in K$ with $\{a,b\}\in E(\comp(P))$. 
     For $0\leq i\leq r$, let $P_i$ be the set of all height $i$ elements of $P$. 
     As $P$ is pure and not an ordinal sum of two induced subposets, $P_i$ is not a singleton for all $i$.
     Clearly, $a,b\notin P_i$ for any $0\leq i\leq r$.
     Assume that $a\in P_i$ and $b\in P_j$ for some $i<j$.
     
     First we show that $P_i \subset K$ and $P_j \subset K$.
     Suppose that there exists an element a $x\in P_i\setminus K$. 
     Let $P_{i,j}$ be the induced subposet of $P$ consisting of height $i$ and height $j$ elements.
     Since $P_{i,j}$ is connected, there exists a sequence $a=a_0,b_0,a_1,b_1,\ldots,b_{s-1},a_s=x$ in $P_{ij}$ such that 
     $a_k\lessdot b_k$ and $a_{k+1}\lessdot b_k$ in $P_{i,j}$ for all $0\leq k\leq s-1$.
     Observe that $a_k \in P_i$ for all $0\leq k \leq s$, and $b_k \in P_j$ for all $0\leq k \leq s-1$. 
     Let $k'$ be the least integer with $0\leq k'\leq s$, such that either $a_{k'}\notin K$ or $a_{k'}\in K$ and $b_{k'}\notin K$.
     If $a_{k'}\notin K$, then $k'\geq 1$ because $a=a_0\in K$.
     Thus, by the choice of $k'$, $ a_{k'-1}, b_{k'-1}\in K$. 
      Also, note that $\{a_{k'},b_{k'-1}\}\in E(\comp(P))$  because $a_{k'}\lessdot b_{k'-1}$ 
     and $\{a_{k'},a_{k'-1}\}\notin E(\comp(P))$ because $a_{k'-1}, a_{k'}\in P_i$. 
 Which is a contradiction to the assumption that $K$ is partitive.
 On the other hand, assume that $a_{k'}\in K$ and $b_{k'}\notin K$. 
 As $b_{k'}\notin K$ and $b\in K$, we get that $b_{k'}\neq b$.
 Note that $\{a_{k'},b_{k'}\}\in E(\comp(P))$ because $a_{k'}\lessdot b_{k'}$ and $\{b, b_{k'}\}\notin E(\comp(P))$ because $b, b_{k'}\in P_j$.
     This again leads to a contradiction to the assumption that $K$ is partitive. 
 Hence, $P_i \subset K$.
     Similarly, one can show that $P_j \subset K$.

     Now we show that if $x\in K\cap P_k$ for some $0\leq k\leq \rank(P)$, then $P_k \subset K$.
     We have already considered the cases when $k=i$ and $k=j$. So we may assume that $k\neq i,j$.
     Since $P$ is pure, there exists a $y\in P_i\subset K$ such that $x$ and $y$ are comparable in $P$.
     So, we can use the above argument to show that $P_k\subset K$.

     Since $K$ is a proper subset of $V(P)$, there exists a $j'$ with $0\leq j'\leq \rank(P)$ such that $K\cap P_{j'} =\varnothing$. Otherwise, if $P_{k}\cap K \neq \varnothing$ for all $k$, then $P_k\subset K$ for all $k$. Thus, $K=V(P) = V(\comp(P))$ which contradicts the assumption $K$ is a proper subset of $V(P)$.
         One can choose $j'$ such that $K\cap P_{j'} =\varnothing$ and either $P_{j'-1}\subset K$ (equivalently, $P_{j'-1}\cap K \neq \varnothing$) or $P_{j'+1}\subset K$ (equivalently, $P_{j'+1}\cap K \neq \varnothing$).
     First assume that $P_{j'-1}\subset K$.
     Since $P$ is not an ordinal sum of two induced subposets, $P_{j'-1,j'}$ is not an ordinal sum of $P_{j'-1}$ and $P_{j'}$.
     Thus, there exists an $x\in P_{j'}$ that does not cover every element of $P_{j'-1}$, say it does not cover $y$. 
     Since $P_{j'-1,j'}$ is pure, $x$ covers some element of $P_{j'-1}$, say $z$.
     Then $\{z,x\} \in E(\comp(P))$ but $\{y,x\} \notin E(\comp(P))$. 
     Hence, $K$ is not partitive.
    Similarly, one can show that if $P_{j'+1}\subset K$, then $K$ is also not partitive.
    Which is a contradiction.
    Hence the proof.
 \end{proof}

\begin{corollary}
	\label{cor:cm_imply_upo}
    Let $G$ be a connected permutation graph. If $G$ is Cohen-Macaulay, then $G^{c}$ is UPO.
\end{corollary}
\begin{proof}
    Let $P$ be a poset such that $G=\cocomp (P)$. Then $\dim P =  2$ as $G$ is connected. As $G$ is Cohen-Macaulay and connected, $P$ is Cohen-Macaulay and not an ordinal sum of two induced subposets. Since $G^{c}=\comp (P)$, the proof follows from the Lemma~\ref{thm:cmpermutation} and Theorem~\ref{thm:upo}.
\end{proof}
	
The converse of Corollary~\ref{cor:cm_imply_upo}, may not be true. The permutation graph in the  Example~\ref{eg:noncm_perm} is connected, and its complement is a UPO, but the graph is not Cohen-Macaulay.

 \subsection*{Acknowledgment}
	The majority of this work was completed during the authors' stay at Chennai Mathematical Institute for the workshop on "Cohen-Macaulay simplicial complexes in graph theory" in July 2023. The authors are thankful to the organizers for their hospitality.
    The last author thanks Ayesha Asloob Qureshi for discussions and her hospitality at Sabancı University, Turkey, in August 2022. This problem originated from those discussions.


\end{document}